\documentclass[10pt]{article}
\newtheorem{theorem}{Theorem}
\newtheorem{lemma}{Lemma}
\newtheorem{remark}{Remark}

\def\Frac#1#2{\frac{\displaystyle{#1}}{\displaystyle{#2}}}
\usepackage{amssymb}
\usepackage{epsfig}

\begin{document}

 \title{Sharp bounds for cumulative distribution functions}

\author{
Javier Segura \thanks{The author acknowledges financial support from Ministerio de Econom\'{\i}a
y Competitividad (project MTM2012-34787)}\\
        Departamento de Matem\'aticas, Estad\'{\i}stica y 
        Computaci\'on,\\
        Universidad de Cantabria, 39005 Santander, Spain.\\
        e-mail: javier.segura@unican.es
}

\date{\ }

\maketitle
\begin{abstract}
Ratios of integrals can be bounded in terms of ratios of integrands
under certain monotonicity conditions. This 
result, related with L'H\^{o}pital's monotone rule, can be used to obtain sharp bounds
for cumulative distribution functions. We consider the case of noncentral cumulative gamma and beta distributions. 
Three different types of sharp bounds for the noncentral gamma 
distributions (also called Marcum functions) are obtained in terms of modified Bessel functions and one additional type of 
function: a second modified Bessel function, two error functions or one incomplete gamma function. 
For the noncentral beta case the bounds are expressed in terms of Kummer functions and one additional Kummer function or an
incomplete beta function. These bounds improve previous results with respect to their range of application and/or its sharpness.
\end{abstract}

2000MSC: 33E20,33B20,26D07,26D15

\section{Introduction}

Given two functions $f_i(x)$ integrable in $[a,b]$ and continuous in $(a,b)$ and
\begin{equation}
\label{def}
F_{i}(x)=\int_a^x f_i(t)dt,\,i=1,2,
\end{equation}
we have that, on account of L'H\^{o}pital's rule, the ratio $F_2(x)/F_1(x)$ has the same limit as $f_2 (x)/f_1(x)$
for $x\rightarrow a^+$ and the same is true for $(F_2(x)-F(b))/(F_1(x)-F(b))$ as $x\rightarrow b^{-}$.
Additionally, as we will see, when $f_2(x)/f_1(x)$ is continuous and monotonic, the ratio 
becomes a bound 
for the ratio of integrals.
This result is closely related to L'H\^{o}pital monotone rule \cite[Theorem 2]{Anderson:2006:MRI} 
(see also \cite[Lemma 2.2]{Anderson:1993:IFQ} and \cite[p. 42]{Cheeger:1982:FPS}).

These bounds on ratios of integrals can be used for bounding cumulative distribution functions 
$G_i(x)=F_i(x)/F_i(b)$ 
when additional information is available, like, for instance, the values
of the difference $G_2(x)-G_1(x)$ or a known explicit expression for one of the functions $G_i (x)$. 
We consider this type of bounds for the noncentral gamma and beta distributions, obtaining different 
families of sharp bounds for different selections of the pair of functions $G_i$.  The noncentral
gamma distributions are also called Marcum functions and have as a particular case (central case)
the incomplete gamma function ratios. On the other hand, the central case for the beta distribution
is the incomplete beta function ratio. 
  
The gamma and beta distributions are classical cumulative distributions appearing in many scientific applications. 
The noncentral gamma distribution is an important function in radar detection and communications
and it is widely used
in statistics and probability theory (see for instance \cite{Segura:2014:MPA} and references cited therein). 
Beta distributions are also important
functions in statistics and probability. The computation of these functions, particularly for the noncentral case,
is difficult. Only recently accurate algorithms for both the lower and upper tail noncentral gamma distributions
(Marcum functions $P$ and $Q$), valid for a large range of parameters, were developed. 
For the noncentral beta case the situation is worse and the available methods
of computation \cite{Chat:1995:ANO} are based on the definition of the lower tail distribution 
in terms of incomplete gamma functions (Eq. (\ref{defi})). 
Probably due to this difficulty, many researchers have been involved in obtaining bounds and  approximations in terms of simpler functions, particularly in the case of the gamma distribution (see, for example, \cite{Baricz:2009:TBF} and references cited therein).

We obtain several types of bounds for the gamma and beta distributions.
For the noncentral gamma  distribution 
we obtain three types of bounds in terms of modified Bessel functions alone, in terms of a modified Bessel function plus two
error functions or in terms of a Bessel function and an incomplete gamma function.
The first type of bounds was recently considered in \cite{Segura:2014:MPA} but using a different approach, while the second type is
 related to some of the bounds in \cite{Baricz:2009:TBF}. 
For the noncentral beta case the bounds are expressed in terms of Kummer functions alone or in terms of an incomplete beta 
function and a Kummer function. In addition to the bounds related to L'H\^opital rule, we will see 
how the recurrence relations satisfied by these functions provide complementary bounds.
 Both for the gamma and beta cases, the bounds we obtain improve previous results with respect to their range of application 
and/or its sharpness. 

We will also prove some 
monotonicity properties and bounds for Kummer functions which are needed in the construction of the bounds for the beta function (see
the appendix). Additionally, we will briefly explain how the bounds are useful for computing the inverse of the distributions (named
quantile functions in statistics).

\section{Bounds related to L'H\^opital's rule}
\label{mainres}

In this section we consider functions defined by integrals. The same results
are valid in general for functions $F_i(x)$  satisfying $F_i(a)=0$ and which are continuous in $[a,b]$ 
and differentiable in $(a,b)$; we are not loosing generality by writing these functions as integrals. We assume that 
$F_i(b)\neq 0$ and we will deal
with normalized integrals $G_i (x)=F_i(x)/F_i(b)$ (so that $G_1(b)=G_2(b)=1$); it is clear how to write
the results for $F_i$ from those for $G_i$.

\subsection{Bounds for ratios of integrals}
\label{ratios}

The following theorem can proved by using Rolle's theorem and Cauchy's mean value theorem.

\begin{theorem}
\label{main}
Let 
$$
G_i(x)=\int_a^x g_i (t)dt,\,i=1,2
$$
with $G_1(b)=G_2(b)=1$, $g_1(x)$ and $g_2(x)$ integrable in $[a,b]$ and continuous in $(a,b)$ and with
$g_1(x)\neq 0$ in $(a,b)$. Let $\bar{G}_i(x)=1-G_i(x)$. 

Then, if $g_2 (x)/g_1(x)$ is strictly
monotonic in $(a,b)$ the following holds:
\begin{enumerate}
\item{There exists} one and only one $x_0$ in $(a,b)$ such that $g_1(x_0)=g_2(x_0)$
\item{}Let $h(x)=G_2(x)/G_1(x)$ and $r(x)=g_2(x)/g_1(x)$. 
The following statements are equivalent in $(a,b)$:
\begin{enumerate}
\item{} $h(x)<1$ ($h(x)>1$)
\item{} $r(x)$ is increasing (decreasing).
\item{} $h(x)<r(x)$ ($h(x)>r(x)$)
\item{} $h(x)$ is increasing (decreasing)
\end{enumerate}
The same holds for $\bar{h}(x)=\bar{G}_2 (x)/\bar{G}_1(x)$ but reversing the inequalities in (a) and (c).
\end{enumerate}
\end{theorem}

{\it Proof.} First we notice that $G_i(x)$ are continuous functions in $[a,b]$ and differentiable in $(a,b)$.

1. Let $D(x)=G_2(x)-G_1(x)$. Because $D(a)=D(b)=0$, by Rolle's theorem there
exists $x_0\in (a,b)$ such that $D'(x_0)=0$, and therefore $g_2 (x_0)=g_1 (x_0)$. And because $g_2(x)/g_1(x)$
is strictly monotonic no other value $x\neq x_0 $ exists such that  $g_2(x)/g_1(x)=1$.

2. We prove the chain of implications $(a)\Rightarrow (b)\Rightarrow (c)\Rightarrow (d)\Rightarrow (a)$ for the case $h(x)<1$; the
case $h(x)>1$ is analogous.

\begin{enumerate}
\item{}$(a)\Rightarrow (b)$: We apply Cauchy's mean value theorem for $G_i(x)$, $i=1,2$ in $[a,x_0]$, with $x_0$ the only value
such that $r(x_0)=1$. Then there exists 
$c\in (a,x_0)$ such that
$$
h(x)=\Frac{G_2(x)}{G_1(x)}=\Frac{G_2 (x)-G_2(a)}{G_1(x)-G_1(a)}=\Frac{g_2(c)}{g_1(c)}=r(c)
$$
and because $h(x)<1$ then $r(c)<1$ with $c<x_0$. And because $r(x_0)=1$ and $r(x)$ is strictly monotonic by hypothesis, then it must be strictly increasing.
\item{}$(b)\Rightarrow (c)$: We consider Cauchy's
mean value theorem in $[a,x]$, $x\le b$: there exist $c\in (a,x)$ such that
\begin{equation}
\label{bo}
h(x)=\Frac{G_2(x)}{G_1(x)}=\Frac{g_2(c)}{g_1(c)}<\Frac{g_2(x)}{g_1(x)}=r(x),
\end{equation}
where the last inequality holds because $c<x$ and $r(x)$ is strictly increasing.
\item{}$(c)\Rightarrow (d)$: Taking the derivative
\begin{equation}
\label{hp}
h'(x)=\Frac{g_1 (x)}{G_1(x)}\left(r(x)-h(x)\right).
\end{equation}
And because $h(x)<r(x)$ then $h(x)$ has positive derivative (because $g_1(x)\neq 0$ for
any $x\in (a,b)$ then $g_1(x)$ and $G_1(x)$ have the same sign). 
\item{}$(d)\Rightarrow (a)$: Because $h(x)$ is strictly increasing in $(a,b)$ and $h(b)=1$ then necessarily $h(x)<1$ in $(a,b)$.
\end{enumerate}

For $\bar{G_i}(x)$ we can proceed analogously. For instance, for the second implication we can apply 
Cauchy's mean value theorem for $G_i(x)$ in $[x,b]$:
$$
\Frac{G_2 (b)-G_2(x)}{G_1 (b)-G_2(x)}=\Frac{\bar{G}_2(x)}{\bar{G}_1(x)}=\Frac{g_2 (c)}{g_1 (c)}>\Frac{g_2 (x)}{g_1 (x)}.
$$
The third implication follows by considering the derivative of $\bar{h}(x)=\bar{G}_2 (x)/\bar{G}_1 (x)$, which is
$$
\bar{h}'(x)=-\Frac{g_1 (x)}{\bar{G}_1(x)}
\left(r(x)-\bar{h}(x)\right).
$$ \hfill $\square$

\begin{remark}
We have proved that $g_2(x)/g_1(x)$ is a bound. In addition, there is the trivial bound $1$ (because $r(x)<1$ or $r(x)>1$). 
Therefore if,
for instance, $h(x)<1$ we have
$$
\Frac{G_2(x)}{G_1 (x)}<\min\left\{1,\Frac{g_2(x)}{g_1 (x)}\right\},\,
\Frac{\bar{G}_2(x)}{\bar{G}_1 (x)}>\max\left\{1,\Frac{g_2(x)}{g_1 (x)}\right\}.
$$

In all cases, the bound $g_2(x)/g_1(x)$ is sharp for $G_2(x)/G_1(x)$ as $x\rightarrow a^+$ and for 
$\bar{G}_2(x)/\bar{G}_1(x)$ as $x\rightarrow b^-$.
The trivial bound $1$ is sharp on the opposite extreme of the interval. The bound $g_2(x)/g_1(x)$ for $G_2(x)/G_1(x)$
($\bar{G}_2(x)/\bar{G}_1(x)$)
is  better than the trivial bound only for $x<x_0$ ($x>x_0$).
\end{remark}

\begin{remark}
That the monotonicity properties of $g_2 (x)/g_1 (x)$ are inherited by $G_2(x)/G_1(x)$ (implication $(b)\Rightarrow (d)$)
is the result named {\it monotone 
L'H\^{o}pital rule} \cite[Theorem 2]{Anderson:2006:MRI}, and it is proved in  \cite[Lemma 2.2]{Anderson:1993:IFQ} similarly as we have proved the two implications $(b)\Rightarrow (c) \Rightarrow (d)$.
\end{remark}

\begin{remark}
An alternative proof for $(a)\Rightarrow (b)\Rightarrow (c) \Rightarrow (d)$ can be given which provides a graphical picture of the situation. 
The proof goes as follows: 

We start from $(a)$ and prove $(b)$ as before. Now, we 
use (\ref{hp}) to conclude that $h'(x)$ has the same sign as $r(x)-h(x)$. With this it is easy to check that
the graph of $h(x)$ is inside the region $(a,b)\times (h(a^+),1)$ and below the graph of $r(x)$. Indeed $h(b^-)<1<r(b^-)$ (because
$h(x)=G_2 (x)/G_1(x)<1$ and $r(b^-)>r(x_0)=1$) and
then, considering (\ref{hp}), $h(x)$ is increasing close to $b$; but it stays
 being increasing in the rest of the interval because, as we move
to the left of the interval, the graph
of $h(x)$ can not cross the graph of the increasing function $r(x)$. If it did at some $x=x_c<x_0$, then $h(x)$ would become decreasing for $x<x_c$, implying that $h(x)>h(x_c)=r(x_c)$. 
Because $x_c>a$ and $r(x)$ is increasing, this would contradict the fact that $h(a^+)=r(a^+)$ on account of L'H\^{o}pital's rule.

\end{remark}

\subsection{Bounds for the functions}
\label{bfu}

The bounds for the ratios can be translated
into bounds for the functions if some of the functions $G_1$ or $G_2$
are explicitly know. We will use this case for obtaining bounds relating the noncentral
and the central distributions; also, for the gamma distribution we will obtain bounds in terms of error functions using this idea.

Another possibility is that there is some relation between $G_1$ and $G_2$ as, for example, 
the difference $G_1(x)-G_2(x)$. This happens both for the gamma and beta distributions, where the recurrence relation they satisfy
provides such a relation.

Considering then that $G_1(x)-G_2(x)$ is known we write:
$$
G_1 (x)=\Frac{G_1(x)-G_2(x)}{1-\Frac{G_2 (x)}{G_1 (x)}}
$$
and considering that the case $h(x)<1$ of theorem \ref{mainres} holds we have 
$G_2(x)/G_1(x)<r(x)=g_2(x)/g_1(x)$, with $g_2(x)/g_1(x)<1$ if $x<x_0$ ($g_2(x_0)/g_1(x_0)=1$). Therefore
\begin{equation}
\label{uno}
G_1(x) < \Frac{G_1(x)-G_2(x)}{1-r(x)}\equiv U_1 (x),\,x<x_0,
\end{equation}
and similarly
\begin{equation}
\label{dos}
G_2< U_2 (x)=r(x) U_1 (x),\,x<x_0.
\end{equation}
These two bounds are sharp as $x\rightarrow a^+$.

In a similar way, we obtain 
\begin{equation}
\label{tres}
\bar{G}_1 (x) < -U_1 (x),\,x>x_0,
\end{equation}
and
\begin{equation}
\label{cuatro}
\bar{G}_2 (x) < -U_2 (x),\,x>x_0,
\end{equation}
and these two bounds are sharp as $x\rightarrow b^-$.

Next we consider the application of Theorem \ref{main} to the gamma 
and beta distributions.

\section{Bounds for gamma distributions}

The Marcum functions (or, equivalently, non-central gamma or chi-squared distributions) are the 
cumulative distribution functions
$$
P_{\mu}(x,y)=\int_0^y g_{\mu}(x,t)dt,\, Q_{\mu}(x,y)=1-P_{\mu} (x,y)
$$
with
$$
g_{\mu}(x,y)=\left(\Frac{y}{x}\right)^{(\mu-1)/2}e^{-x-y}I_{\mu-1}(2\sqrt{xy}).
$$
We have $P_{\mu}(x,0)=1$ and $P_{\mu}(x,+\infty)=1$, as corresponds to a cumulative distribution function in $[0,+\infty)$.

\begin{remark}
We observe that the role of the variable $x$ in section \ref{mainres} is now played by the variable $y$. We use this notation 
in order to be consistent
with previous notation. This notation is used in the rest of this article.
\end{remark}

Integrating by parts we obtain the recursion formulas:
\begin{equation}
\label{recurrgam}
P_{\mu+1}(x,y)=P_{\mu}(x,y)-g_{\mu+1}(x,y),\,Q_{\mu+1}(x,y)=Q_{\mu}(x,y)+g_{\mu+1}(x,y)
\end{equation}

We consider three types of bounds. The first type of bounds were already explored  \cite{Segura:2014:MPA} (but using a different approach); some of the results in that paper are particular cases of these type of bounds.  These bounds involve monotonicity properties
over $\mu$ and a know expression for the difference given by the recurrence relation. The second type of bounds (error function bounds) can be obtained
from the monotonicity over $\mu$ and the particular expression for $\mu=1/2$, $\mu=3/2$. Finally the third type of bounds
(in terms of incomplete gamma functions) are based on the monotonicity with respect to $x$ and the particular value for $x=0$.

\subsection{Bounds for the functions using their recurrence}

The hypotheses of Theorem \ref{ratios} (case $h(x)<1$) hold with the correspondences $G_1\rightarrow P_{\mu}$ 
($\bar{G}_1\rightarrow Q_{\mu}$), $g_1\rightarrow g_{\mu}$, $G_2\rightarrow P_{\mu+1}$ 
($\bar{G}_2\rightarrow Q_{\mu+1}$)
and $g_2\rightarrow g_{\mu+1}$ 
An because $g_{\mu+1}(x,y)/g_{\mu}(x,y)$ is monotonic as a function of $y$ 
\cite[lemma 2]{Segura:2014:MPA}, we have:
$$
\begin{array}{l}
\Frac{P_{\mu+1}(x,y)}{P_{\mu}(x,y)}<\Frac{g_{\mu +1}(x,y)}{g_{\mu}(x,y)}=
\displaystyle\sqrt{\Frac{y}{x}}\Frac{I_{\mu}(2\sqrt{xy})}{I_{\mu -1}(2\sqrt{xy})}\equiv c_{\mu}(x,y)
\end{array}
$$
and
$$
\Frac{Q_{\mu+1}(x,y)}{Q_{\mu}(x,y)}>c_{\mu}(x,y).
$$
The first bound is slightly improved in Theorem 3 of \cite{Segura:2014:MPA} while the second is given in Theorem 6 of that
reference.

From these bounds for the ratios, we can obtain bounds for the function considering the expression given for $G_1-G_2$
which is obtained from the recurrences (\ref{recurrgam}).
Therefore the bounds (\ref{uno}-\ref{cuatro}) hold, with the correspondence 
$G_1(x)-G_2(x)\rightarrow g_{\mu+1}(x,y)$. These are some of the bounds
in section  3.2 of \cite{Segura:2014:MPA}.

From the construction of the bounds (related to L'H\^opital's rule) we deduce they are sharp as $y\rightarrow 0$ for $P_{\mu}$ and
as $y\rightarrow +\infty$ for $Q_{\mu}$. In addition, as explained in \cite{Segura:2014:MPA}, the bound for $P_{\mu}$ is also
sharp as $x\rightarrow +\infty$ and $\mu\rightarrow +\infty$.

Additional bounds can be obtained from the recurrences alone (\ref{recurrgam}), without invoking Theorem \ref{main}. 
Indeed, because the $P$ and the $Q$ functions are
positive $P_{\mu}(x,y)>g_{\mu+1}(x,y)$ and $Q_{\mu}(x,y)>g_{\mu}(x,y)$. We can iterate the recurrence and obtain more sharp bounds. For instance $P_{\mu}(x,y)>\sum_{n=1}^{N}g_{\mu+n}(x,y)$ for any $N$. We refer to \cite{Segura:2014:MPA} for further details.

\subsection{Bounds in terms of error functions}

We take $G_2(y)=\int_0^y g_{\mu+\alpha}(x,t)dt$ and $G_1 (y)=\int_0^y g_{\mu}(x,t)dt$, then
$$
r(y)=\Frac{g_{\mu+\alpha} (x,y)}{g_{\mu}(x,y)}=C(a,x)(\sqrt{y})^\alpha \Frac{I_{\mu-1+\alpha}(2\sqrt{xy})}{I_{\mu-1}(2\sqrt{xy})}.
$$
with $C(a,x)$ not depending on $y$. 
The ratio $r(y)$ is strictly monotonic (increasing) as a function of $y$ for $\mu\ge 1$ as the next lemma shows.

\begin{lemma}
\label{mon1}
The function $t^\alpha \Frac{I_{\mu+\alpha}(t)}{I_{\mu}(t)}$ is increasing as a function of $t$ for $\alpha,\mu>0$.
\end{lemma}

{\it Proof.}
Using the relations in \cite[10.29(i)]{Olver:2010:BES} we obtain:
$$
\left(t^a \Frac{I_{\mu+\alpha}(t)}{I_{\mu}(t)}\right)^{\prime}=t^\alpha \Frac{I_{\mu+\alpha}(t)}{I_{\mu}(t)}\left(\Frac{I_{\mu+\alpha-1}(t)}
{I_{\mu+\alpha}(t)}-
\Frac{I_{\mu-1}(t)}{I_{\mu}(t)}
\right).
$$
And because $I_{\nu -1}(t)/I_{\nu}(t)$ is increasing as a function of $\nu> 0$ 
\cite[Lemma 2]{Segura:2014:MPA} the result is proven. $\square$

Now, using that $r(y)$ is increasing and applying theorem \ref{ratios} we have that

\begin{theorem}
\label{bom}
If $\nu>\mu\ge 1$ then 
$$P_{\nu}(x,y)< \left(\Frac{y}{x}\right)^{\frac{\nu-\mu}{2}} \Frac{I_{\nu-1}(2\sqrt{xy})}{I_{\mu-1}(2\sqrt{xy})}P_{\mu}(x,y)$$
and
$$Q_{\nu}(x,y)> \left(\Frac{y}{x}\right)^{\frac{\nu-\mu}{2}} \Frac{I_{\nu-1}(2\sqrt{xy})}{I_{\mu-1}(2\sqrt{xy})}Q_{\mu}(x,y)$$
\end{theorem}

\begin{remark}
We conjecture that  $I_{\mu-1}(t)/I_{\mu}(t)$, $t>0$ is increasing as a function of $\mu\ge -1/2$ and, therefore, that 
lemma \ref{mon1} holds for $\mu\ge -1/2$, $a>0$ and theorem \ref{bom} holds for $\mu\ge 1/2$. A numerical check shows that 
this is true.
\end{remark}

Taking $\mu=1/2$ in the previous theorem (and taking into account the previous remark) we obtain bounds similar to those in \cite{Baricz:2009:TBF}. 
For instance, and in terms of the variables $a=\sqrt{2x}$, $b=\sqrt{2y}$ (as in \cite{Baricz:2009:TBF}). From the bound for $Q_{\nu}$ in
the previous theorem we have:

\begin{equation}
\label{bo1}
\begin{array}{l}
\tilde{Q}_{\nu}(a,b)>B^{[1]}_{\nu}(a,b)\\
B^{[1]}_{\nu}(a,b)=\Frac{G_{\nu}(a,b)}{\cosh(ab)}\left(\mbox{erfc}\left(\frac{b-a}{\sqrt{2}}\right)+\mbox{erfc}\left(\frac{b+a}{\sqrt{2}}\right)\right)I_{\nu-1}(ab),
\end{array}
\end{equation}
and from the bound for $P_{\nu}$:
\begin{equation}
\label{bo2}
\begin{array}{ll}
\tilde{Q}_{\nu}(a,b)>B^{[2]}_{\nu}(a,b)\\
B^{[2]}_{\nu}(a,b)=1-\Frac{G_{\nu}(a,b)}{\cosh(ab)}\left(\mbox{erf}\left(\frac{b-a}{\sqrt{2}}\right)+\mbox{erf}\left(\frac{b+a}{\sqrt{2}}\right)\right)I_{\nu-1}(ab)
\end{array}
\end{equation}
where $\tilde{Q}_{\nu}(a,b)=Q(a^2/2,b^2/2)$ ($\tilde{Q}_{\nu}(a,b)$ is the function denoted as $Q_{\nu}(a,b)$ in \cite{Baricz:2009:TBF})
$$
G_{\nu}(a,b)=\Frac{a}{2}\displaystyle\sqrt{\Frac{\pi}{2}}\left(\Frac{b}{a}\right)^{\nu}.
$$
The bounds are valid for $\nu>1/2$ and any $a>0$, $b>0$. 

As a marginal note, we notice that the inequality (\ref{bo1}) implies the following relation between error functions and Bessel functions:
$$
\left(\mbox{erfc}\left(\frac{b-a}{\sqrt{2}}\right)+\mbox{erfc}\left(\frac{b+a}{\sqrt{2}}\right)\right)I_{\nu}(ab)<
\displaystyle\sqrt{\Frac{2}{\pi}}\Frac{a^{\nu}}{b^{\nu +1}}(e^{ab}+e^{-ab})
$$
valid for $\nu\ge -1/2$, $a,b>0$.

The bounds (\ref{bo1}) and (\ref{bo2}) can be compared with the results is 
\cite{Baricz:2009:TBF} (eqs. (8) and (17)). For instance, Eq. (8) of the reference is:
$$
\tilde{Q}_{\nu}(a,b)\ge 
\Frac{G_{\nu}(a,b)}{\sinh(ab)}\left(\mbox{erfc}\left(\frac{b-a}{\sqrt{2}}\right)-\mbox{erfc}\left(\frac{b+a}{\sqrt{2}}\right)\right)I_{\nu-1}(ab),
$$
valid for $b\ge a>0$ and $\nu\ge 1$. Which is slightly superior to (\ref{bo1}) when it holds, but of more restricted validity; only
when $a$ is small $a<2$ and $b$ is small (but $b>a$) the advantage becomes noticeable. On the other hand, in the 
comparison between (\ref{bo2}) and Eq. (17) of \cite{Baricz:2009:TBF}, we observe that (\ref{bo2}) is a better bound
except if $a<b$ is quite
close to $b$. 

We can obtain sharper bounds (but of more restricted validity with respect to $\mu$) 
by taking larger $\mu$ in Theorem \ref{bom}. For instance, taking $\mu=3/2$ we obtain bounds valid for $\nu>3/2$
and any $a,b>0$ which are sharper that the bounds of eqs. (8) and (17) of 
\cite{Baricz:2009:TBF}. Using Theorem \ref{bom} together with
the recurrence we have

\begin{equation}
\label{lolo1}
\begin{array}{l}
\tilde{Q}_{\nu}(a,b)> \Frac{a B^{[1]}_{\nu}(a,b)}{b\tanh(ab)}+\left(\Frac{b}{a}\right)^{\nu-1}e^{-(a^2+b^2)/2}I_{\nu-1}(ab)
\end{array}
\end{equation}
and
\begin{equation}
\label{lolo2}
\begin{array}{l}
\tilde{Q}_{\nu}(a,b)> 1-\Frac{a B^{[2]}_{\nu}(a,b)}{b\tanh(ab)}+\left(\Frac{b}{a}\right)^{\nu-1}e^{-(a^2+b^2)/2}I_{\nu-1}(ab)
\end{array}
\end{equation}
valid for $\nu>3/2$.

\subsection{Bounds in terms of incomplete gamma functions}

We take $G_2(y)=\int_0^y g_{\nu}(\rho^2 x,t)dt$ and $G_1 (y)=\int_0^y g_{\nu}(x,t)dt$, with $\rho>0$, $\rho\neq 1$, then
$$
r(y)=\Frac{g_{\nu} (\rho^2 x,y)}{g_{\nu}(x,y)}=C_{\nu}(\rho,x)\Frac{I_{\nu-1}(2\rho\sqrt{xy})}{I_{\nu-1}(2\sqrt{xy})}
$$
where $C_{\nu}(\rho,x)$ does not depend on $y$. Therefore, $r(y)$ has the same monotonicity properties as 
$I_{\mu}(\rho z)/I_{\mu}(z)$ as a function of $z$. Using
\cite[10.29(i)]{Olver:2010:BES} we obtain
$$
\left(\Frac{I_{\mu}(\rho z)}{I_{\mu}(z)}\right)^{\prime}=\Frac{I_{\mu}(\rho z)}{zI_{\mu}(z)}\left(\rho z\Frac{I_{\mu+1}(\rho z)}{I_{\mu}(\rho z)}-z\Frac{I_{\mu+1}(z)}{I_{\mu}(z)}\right),
$$
and because  $xI_{\mu +1}(x)/I_{\mu}(x)$ is increasing as a function of $x$ for
$\mu\ge -1$ (see for instance \cite[Lemma 2]{Segura:2011:BFR}), 
 $r(y)$ is increasing for $\nu\ge 0$ if $\rho>1$ and decreasing if $\rho<1$. Applying theorem \ref{ratios} for $\rho<1$ we obtain
\begin{equation}
P_{\nu}(x,y)<\rho^{\mu-1}e^{-x(1-\rho^2)}\Frac{I_{\nu-1}(2\sqrt{xy})}{I_{\nu-1}(2\rho\sqrt{xy})}P_{\nu}(\rho^2 x,y)
\end{equation}
and
\begin{equation}
Q_{\nu}(x,y)>\rho^{\mu-1}e^{-x(1-\rho^2)}\Frac{I_{\nu-1}(2\sqrt{xy})}{I_{\nu-1}(2\rho\sqrt{xy})}Q_{\nu}(\rho^2 x,y)
\end{equation}
and taking the limit $\rho\rightarrow 0$
\begin{equation}
\label{box}
\begin{array}{l}
P_{\nu}(x,y)<e^{-x}\left(\sqrt{xy}\right)^{1-\nu}I_{\nu-1}(2\sqrt{xy})\gamma_{\nu}(y)\equiv b^{[1]}_\nu (x,y),\\
Q_{\nu}(x,y)>e^{-x}\left(\sqrt{xy}\right)^{1-\nu}I_{\nu-1}(2\sqrt{xy})\Gamma_{\nu}(y)\equiv b^{[2]}_\nu (x,y)
\end{array}
\end{equation}
where 
$$
\gamma_{\nu}(y)=\displaystyle\int_0^y t^{\nu -1}e^{-t}dt,\,
\Gamma_{\nu}(y)=\displaystyle\int_y^{+\infty} t^{\nu -1}e^{-t}dt.
$$
are the lower and upper incomplete gamma functions.

Observe that, as expected, 
$$
\begin{array}{l}
\displaystyle\lim_{x\rightarrow 0^+}b^{[1]}_\nu (x,y)=\Frac{1}{\Gamma (\nu)}\gamma_{\nu}(y)=P_{\nu}(0,y),\\
\\
\displaystyle\lim_{x\rightarrow 0^+}b^{[2]}_\nu (x,y)=\Frac{1}{\Gamma (\nu)}\Gamma_{\nu}(y)=Q_{\nu}(0,y)
\end{array}
$$

Furthermore, because $Q_{\mu}(x,y)\le 1$, the second inequality of (\ref{box}) implies the curious inequality
$$
I_{\nu}(2\sqrt{xy})\Gamma_{\nu +1}(y)<e^x\left(\sqrt{xy}\right)^{\nu},\nu > -1, x,y>0
$$

Bounds in terms of incomplete gamma functions are also given in \cite{Paris:SBF}. Those bounds can be applied for $a=\sqrt{2x}$ or 
$b=\sqrt{2y}$ not much larger than $1$. Our bounds do not have this limitation: The bounds \cite{Paris:SBF} appear to sharper 
for small
$a$ or $b$ but our bounds are sharper for larger parameters.

\section{Bounds for beta distributions}

The cumulative noncentral  beta distribution can be defined as
\begin{equation}
\label{defi}
B_{a,b}(x,y)=e^{-x/2} \displaystyle\sum_{j=0}^{\infty} \Frac{1}{j!}\left(\Frac{x}{2}\right)^j I_y (a+j,b)
\end{equation}
where $I_y (a,b)$ is the central beta distribution
\begin{equation}
\label{cent}
I_y (a,b)=\Frac{1}{B(a,b)}\displaystyle\int_0^y t^{a-1} (1-t)^{b-1} dt,\, B(a,b)=\Frac{\Gamma (a)\Gamma (b)}{\Gamma (a+b)}.
\end{equation}
 A different notation used in the literature (see for instance 
\cite{Chat:1995:ANO}) is $I_x (a,b,\lambda)$, where $B_{a,b}(x,y)=I_y (a,b,x)$. 
The complementary function is defined by $\bar{B}_{a,b}(x,y)=1-B_{a,b}(x,y)$.

Using (\ref{cent}) in (\ref{defi}) we get
\begin{equation}
B_{a,b}(x,y)=\displaystyle\int_0^y g_{a,b}(x,t) dt
\end{equation}
with
\begin{equation}
g_{a,b}(x,y)=\Frac{e^{-x/2}}{B(a,b)}y^{a-1} (1-y)^{b-1} M\left(a+b,a,\Frac{xy}{2}\right)
\end{equation}
and similarly for $\bar{B}_{a,b}(x,y)$ but with the integration going from $y$ to $1$. 

As corresponds to a cumulative distribution function, we have 
$$0=B_{a,b}(x,0)\le B_{a,b}(x,y)\le B_{a,b}(x,1)=1$$ 
for $a,b>0$ and $x,y\ge 0$. Notice that, as for the gamma distribution, the variable $y$ corresponds
to the variable $x$ of section \ref{mainres}; also, as for the gamma distribution, $x$ is called the noncentrality 
parameter. The central case corresponds to $x=0$.

\subsection{Bounds for the functions using their recurrence}

Using the recurrence relations for the central distribution given in \cite[8.17(iv)]{Paris:2010:IGA} together with the
definition (\ref{defi}) it is easy to get analogous relations for the noncentral case. 

Using \cite[8.17.20]{Paris:2010:IGA} we obtain
\begin{equation}
\label{reda}
B_{a,b}(x,y)=B_{a+1,b}(x,y)+\Frac{C_{a,b}(x,y)}{a}M (a+b,a+1,xy/2),
\end{equation}
from \cite[8.17.21]{Paris:2010:IGA}
\begin{equation}
\label{redb}
B_{a,b}(x,y)=B_{a,b+1}(x,y)-\Frac{C_{a,b}(x,y)}{b} M (a+b,a,xy/2),
\end{equation}
and from \cite[8.17.19]{Paris:2010:IGA}
\begin{equation}
\label{redab}
B_{a,b}(x,y)=B_{a-1,b+1}(x,y)-\Frac{C_{a,b}(x,y)}{by} M (a+b,a,xy/2),
\end{equation}
where
\begin{equation}
C_{a,b}(x,y)=e^{-x/2}\Frac{y^a (1-y)^b}{B(a,b)}
\end{equation}

For the complementary function, $\bar{B}_{a,b}(x,y)=1-B_{a,b}(x,y)$ the same recurrences hold but with 
the sign of the inhomogeneous term reversed.
From these recurrences, it becomes clear how the results of section \ref{mainres} can be applied in this case. We give the inequalities
for ratios of $B$-functions; for the $\bar{B}$-function the same holds but with the inequality reversed. We concentrate on
the recurrence relation (\ref{reda}); for the rest of recurrences, a similar analysis is possible.

From (\ref{reda}) we observe that the role of $G_1$ and $G_2$ 
in section \ref{mainres} are played by $B_{a,b}$ and $B_{a+1,b}$.
Therefore
\begin{equation}
\label{cotasco}
\begin{array}{l}
\Frac{B_{a+1,b}(x,y)}{B_{a,b}(x,y)}<\Frac{g_{a+1,b}(x,y)}{g_{a,b}(x,y)}=\Frac{a+b}{a}y
\Frac{M(a+b+1,a+1,z)}{M(a+b,a,z)}
\end{array}
\end{equation}
where $z=xy/2$, provided that the function ratio on the right-hand side is monotonic, which is obviously true for the central case $x=0$.
The proof of monotonicity can be obtained by using similar ideas 
as those used in
\cite{Segura:2011:BFR} for the ratios of modified Bessel functions. In the Appendix 
we prove that $z M(a+b+1,a+1,z)/M(a+b,a,z)$ is increasing; in addition, using the bound given in this theorem, we obtain
the simpler bound
$$
\Frac{B_{a+1,b}(x,y)}{B_{a,b}(x,y)}<
\Frac{1}{x}\left(z+1-a+\sqrt{(z+1-a)^2+4(a+b)z}\right),\,z=xy/2
$$

Now we obtain bounds of incomplete beta functions from the bounds for the ratios as described in section \ref{bfu}.
Using (\ref{reda}) and (\ref{cotasco}) we have the bound for 
$G_1$
\begin{equation}
\label{1b}
B_{a,b}(x,y)<\Frac{C_{a,b}(x,y) M (a+b,a+1,xy/2)}{a\left(1-\Frac{g_{a+1,b}}{g_{a,b}}\right)}
\end{equation}
and for $G_2$
\begin{equation}
\label{2b}
B_{a+1,b}(x,y)<\Frac{C_{a,b}(x,y) M (a+b,a+1,xy/2)}{a\left(\Frac{g_{a,b}}{g_{a+1,b}}-1\right)}.
\end{equation}

Let us particularize these bounds for the central case ($x=0$). The inequality (\ref{1b}) gives
\begin{equation}
\label{uu}
B_{a,b}(0,y)<\Frac{y^a (1-y)^b}{B(a,b)(a-(a+b)y)}, y<\Frac{a}{a+b}
\end{equation}
and (\ref{2b})
\begin{equation}
\label{ee}
B_{a,b}(0,y)<\Frac{y^a (1-y)^b}{B(a,b) (a-1-(a-1+b)y)}, y<\Frac{a-1}{a-1+b}
\end{equation}
which is essentially the same bound.

For $\bar{B}_{a,b}(0,y)$, we have similarly to (\ref{uu}) 
\begin{equation}
\label{rela}
\bar{B}_{a,b}(0,y)<\Frac{y^a (1-y)^b}{B(a,b)((a+b)y-a)}, y>\Frac{a}{a+b}
\end{equation}
and the bound related with (\ref{ee}).

The bound (\ref{1b}) compares favourably with the upper bound in \cite[Thm. 5]{Cerone:2007:SFA}. It is simpler and sharper except
in some cases when $a$ is small ($a<3$). For large values of $a$ and/or $b$ our bound is several orders of magnitude sharper.

\subsubsection{Bounds using only the recurrence}

The recurrence relations, without invoking Theorem \ref{main}, provide sharp bounds themselves. 
For instance, considering (\ref{reda}) we conclude that
$$B_{a,b}(x,y)=\Frac{C_{a,b}(x,y)}{a}M(a+b,a+1,xy/2)$$
and for the complementary function we have
$$\bar{B}_{a,b}(x,y)>\Frac{C_{a-1,b}}{a-1} M (a+b-1,a,xy/2).$$

For the case of $B_{a,b}(x,y)$ we can construct a convergent sequence of lower bounds by applying the recurrences consecutively. In that 
way, after applying $N$ times the recurrence we have
$$
B_{a,b}(x,y)=B_{a+N+1,b}+\Frac{C_{a,b}(x,y)}{a} \displaystyle\sum_{j=0}^N\Frac{(a+b)_j}{(a+1)_j}y^jM(a+b+j,a+j+1,xy/2)
$$
which gives the bound
\begin{equation}
\label{ubob}
B_{a,b}(x,y)> \Frac{C_{a,b}(x,y)}{a}\displaystyle\sum_{j=0}^N\Frac{(a+b)_j}{(a+1)_j}y^jM(a+b+j,a+j+1,xy/2)
\end{equation}

\subsubsection{An application: quantile function estimation for the central beta case}

The previous bounds can be used as estimation for the cumulative distribution, particularly when they are sharper (at the tails). 
But they can also be used for computing estimations of the inverse distribution or quantile function. For instance,
the bounds for the central gamma distribution were used in some cases as starting values for the inversion method described in 
\cite{Gil:2015:GAP}.
Similarly, the bounds for the central beta distribution are useful as starting values for inverting the equation $B_{a,b}(0,y)=\beta$,
with $\beta$ small.

Considering the estimation (\ref{uu}) and because $B_{a,b}(0,y)$ is increasing as a function of $y$ and 
$$f(y)=\Frac{y^a (1-y)^b}{B(a,b)(a-(a+b)y)}$$ 
is also increasing as a function of $y$ if $y<a/(a+b)$, then
the solution of $f(y)=\beta$ in $(0,a/(a+b))$ gives a value $y_{l}$ such that $B_{a,b}(0,y_l)<\beta$. Therefore
the value $y_\beta$ such that $B_{a,b}(0,y_\beta)=\beta$ is such that $y_l<y_{\beta}$ and it should be close to $y_l$ if
$\beta$ is sufficient small. Depending on the values of $a$ and $b$, the lower bound $y_l$ guarantees convergence of the fourth order 
Schwartzian-Newton method \cite{Segura:2015:TSN} to $y_{\beta}$; and when it does not guarantee convergence, an upper bound will do.

An upper bound for $y_{\beta}$ can be obtained from the lower bound for $B_{a,b}(x,y)$ (\ref{ubob}). For instance, considering $N=1$ (and $x=0$)
we have 
$$
B_{a,b}(0,y)>\Frac{y^a (1-y)^b}{aB(a,b)}\left(1+\Frac{a+b}{a+1}y\right)=g(y)
$$
and the solution of $g(y)=\beta$ gives an upper bound for $y_\beta$, $y_u>y_{\beta}$. Therefore $y_{\beta}\in (y_l,y_u)$.

The equations $f(y)=\beta$ or $g(y)=\alpha$ can be solved by any numerical method for nonlinear equations. A possibility with good 
convergence properties is to use the fixed point method obtained by isolating $y$ as follows
$$
y=h(y)=(\beta B(a,b)(a-(a+b)y)(1-y)^{-b})^{1/a};
$$
we can iterate $y_{n+1}=h(y_n)$ starting from $y_0=0$ in order to solve numerically the
equation $f(y)=\beta$. We can proceed similarly for $g(y)=\beta$ taking now
$$
h(y)=\left(1+\Frac{a+b}{a+1}y\right)^{-1/a}(a\beta B(a,b)(1-y)^{-b})^{1/a};
$$

These approximations are used as starting values for the the inversion of the cumulative beta distribution for small values
of $\beta$ \cite{Gil:2015:EAF}. One of them (either $y_l$ or $y_u$) guarantees convergence of the fourth order 
Schwartzian-Newton method to $y_{\beta}$.

\subsection{Bounds  in terms of the central distribution}

The roles of $G_1$ and $G_2$ are now played by $g_{a,b}(\rho x,y)$ and $g_{a,b}(x,y)$. We have
$$
\Frac{g_{a,b}(\rho x,y)}{g_{a,b}(x,y)}=C(a,b,x)\Frac{M(a+b,a,\rho z)}{M(a+b,a,z)},\,z=xy/2
$$
and $C(a,b,x)$ does not depend on $y$. Therefore the monotonicity properties of $r(y)=g_{a,b}(\rho x,y)/g_{a,b}(x,y)$ as a function of $y$ are
the same as the properties as a function of $z$ of
 $q(z)=M(\alpha,\beta,\rho z)/M(\alpha,\beta,z)$, $\alpha=a+b$, $\beta=a$. Now, we the aid of \cite[13.3]{Daalhuis:2010:CHF} we obtain
$$
q'(z)=\Frac{\alpha}{z}\Frac{M(\alpha,\beta,\rho z)}{M(\alpha,\beta,z)}
\left[\Frac{M(\alpha+1,\beta,\rho z)}{M(\alpha,\beta,\rho z)}-\Frac{M(\alpha+1,\beta,z)}{M(\alpha,\beta,z)}\right].
$$

In the appendix we will prove that $\Frac{M(\alpha+1,\beta,z)}{M(\alpha,\beta,z)}$ 
is increasing as a function of $z$ for $\alpha,\beta>0$. Therefore
$r(y)$ is decreasing as a function of $y$ for $\rho<1$ and we have

$$
\Frac{B_{a,b}(\rho x,y)}{B_{a,b}(x,y)}>\Frac{g_{a,b}(\rho x,y)}{g_{a,b}(x,y)}=e^{-x(\rho -1)/2}
\Frac{M(a+b,a,\rho xy/2)}{M(a+b,a,xy/2)}
$$

An particularizing for $\rho =0$, we have
$$
B_{a,b}(x,y)<e^{-x/2}M(a+b,a,xy/2)B_{a,b}(0,y)
$$
Similarly
$$
\bar{B}_{a,b}(x,y)>e^{-x/2}{M(a+b,a,xy/2)}\bar{B}_{a,b}(0,y),
$$
which implies that
$$
\bar{B}_{a,b}(0,y)M(a+b,a,xy/2)<e^{x/2}
$$
for $a,b,x>0$, $y\in (0,1)$.

\section*{Appendix}

In order to obtain a bound for the beta distributions, we used some monotonicity properties
for the ratios of Kummer functions that we next prove. 

First, we used the property that
$zM(a+1,b+1,z)/M(a,b,z)$ is increasing as a function of $z$. Next we prove this
result and we supplement it with additional information.

\begin{theorem}
The function
$$h(z)=\Frac{M(a+1,b+1,z)}{M(a,b,z)},$$
is decreasing as a function of $z>0$ if $a>b>0$, constant if $a=b>0$ and
increasing if $0<a<b$, while $zh(z)$ is increasing for all $a,b>0$.

If $a>b>0$ the following holds
$$
D(a,b,z)<\Frac{a}{b}h(z)< D(a,b-1,z)
$$
where
$$
D(a,b,z)=\Frac{1}{2z}\left(z-b+\sqrt{(z-b)^2+4az}\right).
$$
The upper bound also holds for $0<a\le b$, but a sharper bound is given by
$$
\Frac{a}{b}h(z)< D(a,b,z)
$$
\end{theorem}

{\it Proof.} We start with the function $zh(z)=\frac{b}{a}H(z)$ with
$$
H(z)=z\Frac{M'(a,b,z)}{M(a,b,z)}
$$
taking the derivative and using the Kummer equation $zM''(a,b,z)+(b-z)M'(a,b,z)-aM(a,b,z)=0$ we find
$$
H'(z)=a+\left(1+\Frac{1-b}{z}\right)H(z)-\Frac{1}{z}H(z)^2.
$$
We compute the characteristic roots $\lambda(z)$ (such that if $H(z)=\lambda (z)$ then 
$H'(z)=0$) by solving $a+\left(1+\Frac{1-b}{z}\right)\lambda (z)-\Frac{1}{z}\lambda(z)^2=0$.
Because $a>0$ one solution is positive and the other negative. The interesting solution
is the positive, which is $\lambda_+ (z)=zD(a,b-1,z)$.

Because $a,b>0$, $H(z)>0$. We will have $H'(z)>0$ when $0<H(z)<\lambda_+ (z)$ and $H'(z)<0$ when $H(z)>\lambda_+ (z)$,
but the second possibility can be easily ruled out.

We have $\lambda_+ (0^+)=0^+$ and $\lambda_+'(z)>0$, $z>0$. On the other hand, $H(0^+)=0^+$ and $H'(0^+)=1$, which
implies that $H(0^+)<\lambda_+(0)$, but then necessarily $H(z)<\lambda_+(z)$ and $H'(z)>0$ for all $z>0$. Indeed, the
graph of $H(z)$ can not intersect the graph of the characteristic root $\lambda_+ (z)$ because $\lambda_+ (z)$
is increasing and $H(0^+)<\lambda_+(0)$. 

This proves the monotonicity properties of $zh(z)$ and the upper bound of the theorem.

Now we consider $G(z)=\Frac{M'(a,b,z)}{M(a,b,z)}=\Frac{a}{b}h(z)$. Taking the derivative and using Kummer's equation
$$
G'(z)=\Frac{a}{z}+\left(1-\Frac{b}{z}\right)G(z)-G(z)^2.
$$
For $a>0$ there is only one positive characteristic root, which is $\lambda_+ (z)=D(a,b,z)$, which is decreasing if $a>b>0$, equal
to $1$ if $a=b$ and increasing if $0<a<b$; in addition $\lambda_+ (0^+)=a/b$. On the other hand $G(0^+)=a/b$ with
$G'(0^+)<0$ if $a>b$, $G'(z)=0$ if $a=b$ and $G'(0^+)>0$ if $a<b$. With this, and using similar qualitative arguments
as before, the monotonicity properties for $h(z)$ and bound is proved (which gives the lower bound in the theorem).
$\square$

For obtaining the bounds for small $x$, we used monotonicity of the ratio $M(a+1,b,x)/M(a,b,x)$. We prove this result and
an associated bound:

\begin{theorem}
The function
$$\Frac{M(a+1,b,z)}{M(a,b,z)}$$
is increasing as a function of $z$ for $a,b,z>0$ and the following bound holds for that range of the variables:
$$
\Frac{M(a+1,b,z)}{M(a,b,z)}<\Frac{1}{2a}\left(z+1-b+2a+\sqrt{(z+1-b)^2+4za}\right)
$$
\end{theorem}

{\it Proof.}
Given $H(z)=\Frac{M(a+1,b,z)}{M(a,b,z)}$, taking the derivative and using rules in \cite[13.3]{Daalhuis:2010:CHF}
$$
H'(z)=\Frac{1}{z}\left(a+1-b+(z+2a+1-b)H(z)-aH(z)^2\right).
$$
The characteristic roots of this Riccati equation are the roots of $a+1-b+(z+2a+1-b)\lambda(z)-a\lambda (z)^2=0$, which are
$$
\lambda_{\pm}(z)=\Frac{1}{2a}\left(z+1-b+2a\pm \sqrt{(z+1-b)^2+4za}\right).
$$ 
For $z>0$ we have $\lambda_+ (z)>\lambda_{-}(z)$ and $\lambda_+^{\prime}(z)>0$ if $a>0$.
These properties, together with the fact that $H(z)=1+z/b+{\cal O}(z^2)$ as $z\rightarrow 0$ and 
$H(z)\sim\lambda_{+}(z)\sim z/a$ as
$z\rightarrow +\infty$ (see \cite[13.7.1]{Daalhuis:2010:CHF}) are enough to prove the theorem.

Indeed, because $H(z)$ is increasing for sufficiently small $z$  (because $H(z)=1+z/b+{\cal O}(z^2)$)
then $\lambda_{-}(z)<H(z)<\lambda_{+}$ for $z$ small (because $H'(z)<0$ if $H(z)>\lambda_+ (z)$ or $H(z)>\lambda_+ (z)$); the
same is true for large $z$ because $H(z)\sim z/a$ as
$z\rightarrow +\infty$. But then necessarily $\lambda_{-}(z)<H(z)<\lambda_{+}(z)$ for all $z>0$ for the following reasons:

1. 
Let $z_c$ be the smallest value of $z$ such that $H(z_c)=\lambda_+(z_c)$; we prove that such value can not exist. Indeed, if such $z_c$ existed we would have
$H'(z_c)=0$, $H(z_c)=\lambda_+(z_c)$, $H(z)<\lambda_+ (z)$ for $z<z_c$ and with $\lambda_+^{\prime}(z)>0$ for $z>0$, which is not 
possible. Therefore $H(z)<\lambda_+ (z)$ for all $z>0$.

2. For large enough $z$ we have  $H(z)>\lambda_{-} (z)$ but then this holds for all $z>0$. Observe that $\lambda_{-}(z)$ is
monotonic (increasing if $a<b-1$, decreasing of $a>b-1$, constant if $a=b-1$). Then, if $\lambda_{-}(z)$ is increasing or constant
there can not exist $z_c$ such that $H(z_c)=\lambda_{-}(z_c)$ using the same arguments as in the item 1 of this proof. And
if $\lambda_{-}(z)$ is decreasing such number $z_c$ can not exist either: if such a number existed then 
$H(z_c^{-})<\lambda_{-}(z_c^-)$
and $H'(x_c^{-})<0$ but then this implies, because $\lambda_- (x)$ is decreasing, 
that $H'(z)<0$ for $0<z<z_c$, which is not true ($H(z)$ is increasing close to $z=0$).

Then, because $\lambda_{-}(z)<H(z)<\lambda_{+}(z)$ for all $z>0$ we have $H'(z)>0$ for all $z>0$
$\square$

\bibliographystyle{elsart-num-sort}
\bibliography{beta_arxiv}

\end{document}